\documentclass[num-refs]{wiley-article}

\usepackage[frozencache=true,cachedir=minted-cache]{minted} 
\usepackage{siunitx}
\setminted{baselinestretch=0.7}
\newenvironment{myquote}%
  {\list{}{\leftmargin=0.1in\rightmargin=0.1in\color{darkgray}\itshape}\item[]}%
  {\endlist}
\papertype{Focus Article}

\title{Understanding Automatic Differentiation Pitfalls}

\author[1]{Jan~H{\"u}ckelheim}
\author[2]{\hspace{-0.6em}Harshitha~Menon}
\author[3]{\hspace{-0.6em}William~Moses}
\author[4]{\hspace{-0.6em}Bruce~Christianson}
\author[1]{\hspace{-0.6em}Paul~Hovland}
\author[5]{\hspace{-0.6em}Laurent~Hasco{\"e}t}

\affil[1]{Argonne National Laboratory, 9700 S. Cass Ave., Lemont, IL 60439, USA}
\affil[2]{Lawrence Livermore National Laboratory, 7000 East Ave., Livermore, CA 94551, USA}
\affil[3]{MIT CSAIL, 32 Vassar St, Cambridge, MA 02139, USA}
\affil[4]{University of Hertfordshire, College Lane, Hatfield, AL10 9AB, UK}
\affil[5]{Inria, 2004 Rte des Lucioles, 06902 Valbonne, France}

\corraddress{Jan H\"uckelheim, Argonne National Laboratory, Lemont, IL 60439, USA}
\corremail{jhueckelheim@anl.gov}

\fundinginfo{See acknowledgment section.}

\runningauthor{J. H\"uckelheim, H. Menon, W. Moses, B. Christianson, P. Hovland, L. Hasco\"et -- Understanding AD Pitfalls}

\begin{document}

\begin{frontmatter}
\maketitle

\begin{abstract}
Automatic differentiation, also known as backpropagation, AD, autodiff, or algorithmic differentiation, is a popular technique for computing derivatives of computer programs accurately and efficiently.

Sometimes, however, the derivatives computed by AD could be interpreted as incorrect. These pitfalls occur systematically across tools and approaches. In this paper we broadly categorize problematic usages of AD and illustrate each category with examples such as chaos, time-averaged oscillations, discretizations, fixed-point loops, lookup tables, and linear solvers. We also review debugging techniques and their effectiveness in these situations.

With this article we hope to help readers avoid unexpected behavior, detect problems more easily when they occur, and have more realistic expectations from AD tools.

\end{abstract}
\end{frontmatter}

\section{Introduction}
\label{sec:intro}
Automatic differentiation (AD) is steadily becoming more popular in machine learning, scientific computing, engineering, and many other fields as a tool to compute derivatives efficiently and accurately.
While the benefits are widely popularized, users are not always aware that AD can produce surprising results when applied to certain functions. Some of the most problematic failure modes are inherent to common AD approaches and systematically lead to confusing (or, by some measure, incorrect) derivatives across multiple tools or input languages. This makes debugging difficult without awareness and in-depth understanding by AD users.

This article presents a new way to categorize and understand known AD pitfalls. It also gives an overview of debugging techniques that are commonly used by practitioners and shows which techniques are able to identify which types of problems. The goal is to enable readers of this article to use AD more wisely, avoid known pitfalls, and be less surprised when they discover new ones in the future.

\section{What do we expect from AD?}
\label{sec:whattoexpect}
Before discussing pitfalls, we review quotes from influential prior work stating the desired outcome of using AD.

\begin{myquote}
    Frequently we have a program that calculates numerical values for a function, and we would like to obtain accurate values for derivatives of the function as well.~\cite{Griewank:2008:EDP:1455489}

    [AD] describes the mathematical theory of how a computer program can be differentiated.~\cite{codipack}

    [AD] ... is a family of techniques ... for efficiently and accurately evaluating derivatives of numeric functions expressed as computer programs.~\cite{baydinsurvey}

    AD differentiates what you implement ...
    Which occasionally differs from what you think you implement!~\cite{naumann2012art}
\end{myquote}
A reader might wonder what function is being differentiated or how we define the derivatives of a computer program. This question indeed lies at the heart of many pitfalls discussed in this paper.
Let us consider a mathematical function
$y \leftarrow f(x)$,
where $x, y$ are real or complex scalar values or vectors of arbitrary size. One might wish to compute $\nabla(f)$, the derivative of $y$ with respect to $x$, which---depending on the shape and size of $x$ and $y$---could be a scalar, column or row vector, or, more generally, a Jacobian matrix. Instead of computing an entire Jacobian matrix, one might wish to compute projections of that matrix (popularly implemented as the forward mode of automatic differentiation) or of its transpose (implemented as reverse mode or backpropagation). AD can also compute Hessians or other higher-order derivatives or their projections. The pitfalls described in this work apply to all of these AD usages; and we refer to some of the many articles, books, and reviews that discuss AD modes and their relationship~\cite{racc,BARTHOLOMEWBIGGS2000171,griewank_2003,Griewank:2008:EDP:1455489,naumann2012art,hitchhiker,baydinsurvey,margossian2019review,wileyad,yolo}.

As alluded to by the above quotes, instead of the \emph{true} mathematical function $f$, a computer program implements an approximation $Y\leftarrow F(X)$, where $F$, $X$, and $Y$ differ from the true $f$, $x$, and $y$ for reasons that include floating-point errors and often many other approximations. AD may not encounter the true function $f$ and can thus not be expected to compute its derivatives. We can at best hope to obtain derivatives of $F$ or, in other words, the ``program" or ``what you implement." Even if the true and implemented functions produce very similar results for all inputs, their derivatives might be completely different for some or all inputs.

Moreover, AD operates on some level of abstraction, which is determined by the user within the constraints placed by a given AD tool. The level of abstraction influences the exact function $F$ that we consider the program to implement. For example, a program could consist of a few calls to linear algebra kernels. These kernels might internally use iterative solvers or approximations using polynomials and lookup tables. All of this is executed by using floating-point numbers instead of real numbers. Each level of abstraction in this example represents a different function; and while these functions usually have reasonably similar outputs by design, their derivatives might be vastly different.
With this situation in mind, we can attempt to define the expected output of AD as follows:
\begin{myquote}
     AD operates on a chosen level of abstraction, with the assumption that all operations and their derivatives can be computed to sufficient accuracy when expressed at that level of abstraction. 
     Under these assumptions, AD computes the derivative of a function $F$ given as the composition of operations encountered during the evaluation of a particular branch of a program above the given abstraction level.
\end{myquote}
This definition allows us to more clearly categorize and discuss pitfalls in Section~\ref{sec:sources}.
Before doing so, we  note some problems that are outside the scope of this article. The level of abstraction at which AD is applied  influences not only the results but also the run time, memory footprint, user convenience, and tool development effort. On these grounds, some authors argue in favor of AD on high abstraction levels~\cite{doi:10.1137/120873558} while others argue for AD on low abstraction levels~\cite{moses2020instead}. Some languages have operators whose derivatives cannot be efficiently expressed within the same language, and some language constructs (such as asynchronous parallelism, operations involving void pointers, aliasing) are difficult to statically analyze and differentiate~\cite{HascoetUtke2016}. Some operations could have  simple and cheap-to-compute derivatives, but tools fail to optimize expressions to take advantage of this fact. As important as these aspects are, we nonetheless focus on the semantics of AD in this work and ignore efficiency questions.

We also note that AD tools have bugs and limitations just like any other nontrivial software, some of which cause tools to silently produce wrong results. Limitations regarding the supported input programs are sometimes poorly documented or in flux, and different tools and approaches represent different trade-offs regarding performance or language coverage, as reviewed, for example, in~\cite{margossian2019review}. However, this article focuses on problems that occur even with a hypothetical bug-free AD tool applied to a function implemented in a supported input language.

\section{Categories of AD pitfalls}
\label{sec:sources}

In this section we illustrate each pitfall category with small examples and refer to related work for more details. The pitfalls are ordered from high to low level:
We begin with problems caused by the mathematical function, then discuss the role of the implementation at the chosen level of abstraction, and conclude with operator precision. 


\subsection{Pitfall I: The function has unexpected derivatives}
\label{sec:funcdervs}

The derivatives of the true function $f$ do not always exist or are not always useful.
For example, a program may model discontinuous physical processes such as phase changes from liquid to solid, or shocks in fluid dynamics. 
Sometimes the discontinuities are introduced when modeling smooth processes, for example when quantities are filtered to avoid nonphysical behavior or numerical instabilities~\cite{sweby1984high}. Discretized models are often designed to be good approximations of a continuous system, but their derivatives might be vastly different~\cite{Hager2000RungeKuttaMI,collis_adjoints}. In some cases  the derivatives of the discretized system are useful, and in others a discretization of a differentiated system is preferable~\cite{nadarajah2000comparison}. The computational fluid dynamics community--being an early user of AD and adjoint methods--has studied these issues in  detail~\cite{appel1997difficulties,Fikl_2016,doi:10.1137/080727464,derivativesshocks}, but practitioners in other domains will likely discover similar problems when they adopt AD.

\begin{figure}
    \centering
    \includegraphics[scale=0.7]{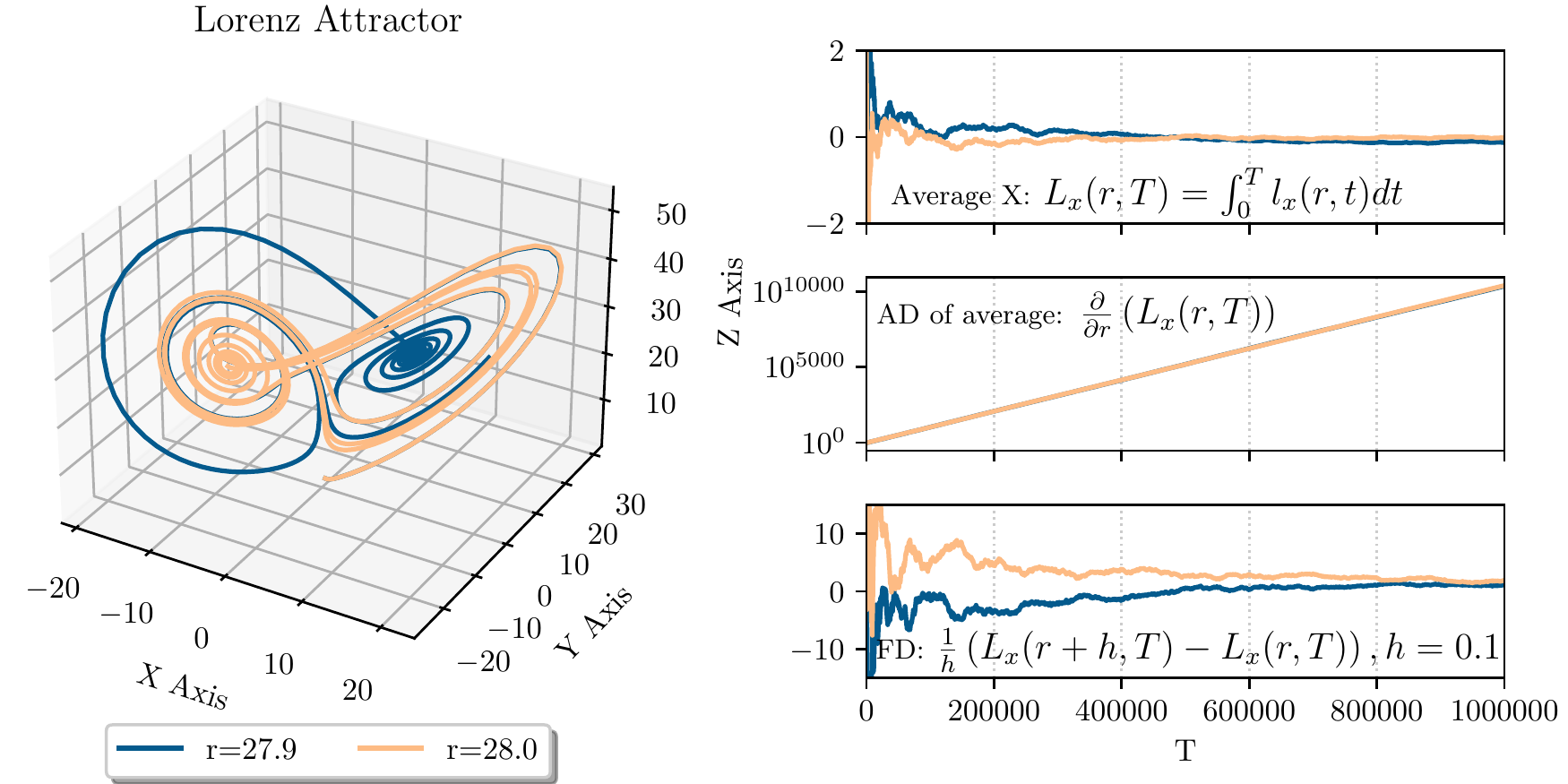}
    \caption{\textbf{Left:} Chaotic functions, such as the Lorenz attractor $l(r, t)$, have ill-behaved derivatives. \textbf{Top right:} The time-averaged behavior converges and becomes a smooth function in $r$ in the limit of $T\rightarrow\infty$. \textbf{Middle right:} For any finite time window, however, the time-averaged function is noisy in $r$, and its derivatives diverge exponentially with growing $T$. \textbf{Bottom right:} Finite differences (FD) do not resolve the noise in the time-averaged function and produce derivatives that are more intuitive (even if less accurate than those produced by AD).}
    \label{fig:lorenz}
\end{figure}

Differentiable functions with well-understood large-scale behavior may have localized behavior that causes unintuitive, noisy, or misleading derivatives. Previous work has found that finite differences  can in such cases be more reliable~\cite{MORE2014268,doi:10.1080/10556780500151984}. One example is the computation of time averages or other statistical properties of dynamic or chaotic processes, which can have exponentially diverging derivatives~\cite{wang2013forward} for finite-size time windows. We illustrate such a case in Figure~\ref{fig:lorenz}.
From a user's perspective, it can be difficult to determine whether  a program (or a small part of it) implements a chaotic function, which may appear as an unremarkable composition of differentiable operations. Consequently, a user may  become aware of small-scale chaotic effects only when applying AD.
Previous work proposed methods to obtain meaningful derivatives in the presence of chaos~\cite{wang2013forward,wang2014least,blonigan2017adjoint}, which go beyond the standard definition of AD and need to be used judiciously because of their high cost.

\begin{figure}[t]
    \centering
    \includegraphics[scale=0.7]{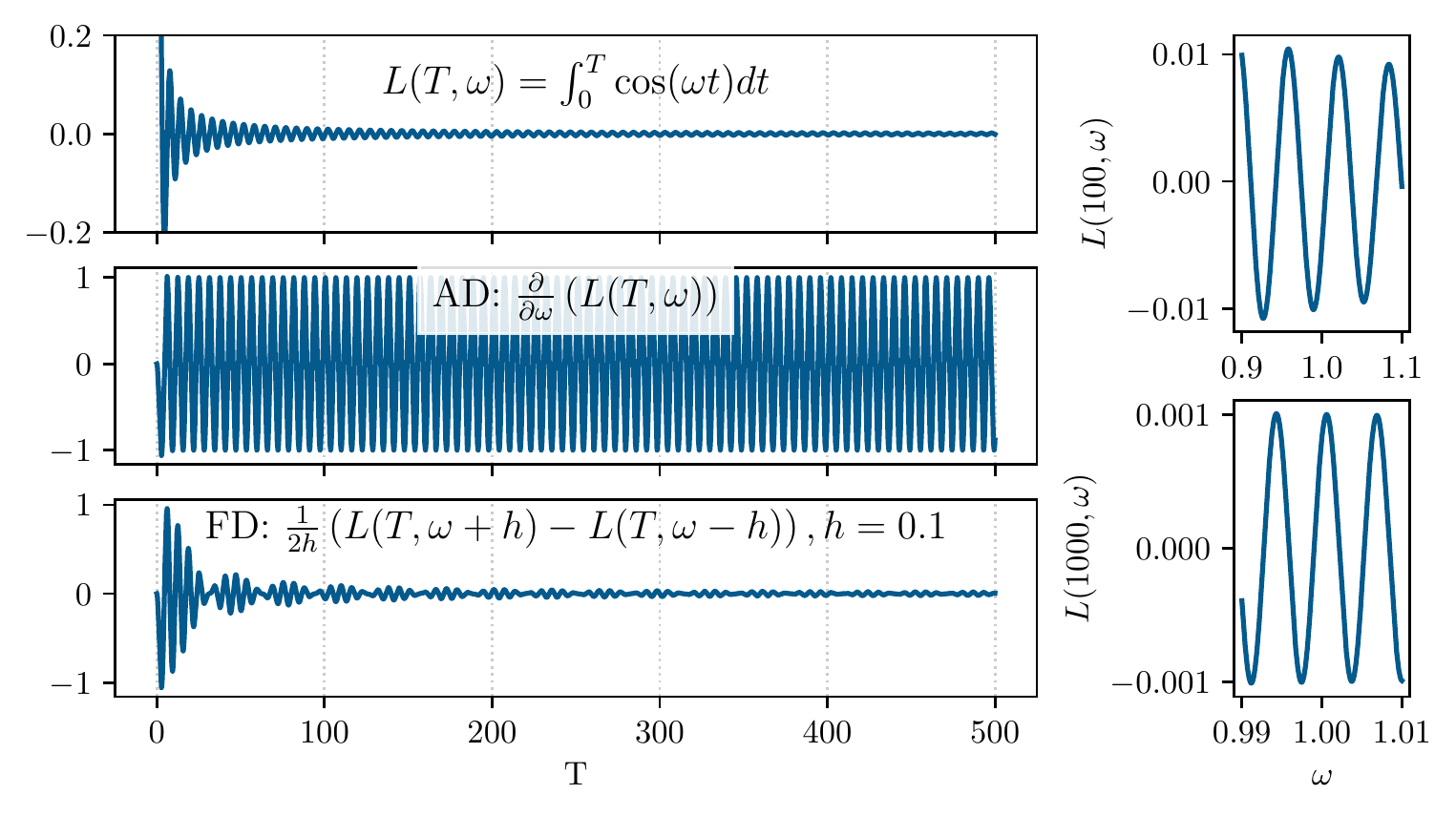}
    \caption{A nonchaotic function with unintuitive derivatives. \textbf{Top left:} The time average converges to zero, regardless of the frequency. One might thus expect the derivative with respect to the frequency parameter to converge to zero. \textbf{Middle left:} The derivative, however, continues to oscillate between $-1$ and $1$. \textbf{Bottom left:} Finite differences (FD) behave more intuitively and converge to zero. \textbf{Right:} As $T$ increases, the amplitude of $L$ decreases but becomes more sensitive to changes in $\omega$. Note the differences in axis scale for the top and bottom right plots.}
    \label{fig:cosine}
\end{figure}

A poorly chosen objective function can itself cause misleading derivatives. Consider a system whose behavior in time is modeled as a cosine function with frequency $\omega$. A user may attempt to compute the time-averaged state by integrating over a long, but finite, time window $[0, T]$. While the integral converges to the same value with growing $T$ regardless of $\omega$ (which might intuitively indicate a $0$-derivative with respect to $\omega$), the derivative never converges and is a plausible-looking but arbitrary value between $-1$ and $1$ for any finite $T$.
Confusingly, finite differences converge to the intuitive, but incorrect, value of $0$ for any fixed step size $h$, because of a failure to resolve the small-scale oscillations. Figure~\ref{fig:cosine} shows this in detail.
The example may appear contrived but represents a challenge in real applications, such as small oscillations occurring during computational fluid dynamics simulations that have been found to cause misleading derivatives for time-averaged objectives like drag or fuel consumption in aerospace engineering. Instead of computing averages with a sharp cutoff at some start and end point, modified objective function formulations with smoothly differentiable time windows~\cite{KRAKOS20123228} can be used to obtain more intuitive derivatives.

Undesirable derivatives are not limited to simulations of physical systems. For example, machine learning models may have ``exploding'' or ``vanishing'' gradients that grow or shrink exponentially with the number of layers and can be avoided by carefully choosing the hidden layers~\cite{vanishing_exploding,vanishing}.
In summary, we observe that useful functions, even differentiable ones, do not always have useful derivatives. The successful use of AD can require a deep understanding of the actual problem being modeled, as well as models and objective functions designed with differentiation in mind.

\subsection{Pitfall II: The chosen abstraction has unexpected derivatives}
\label{sec:abstractiondervs}

Most programs are written at a relatively high level of abstraction and are successively transformed, or \emph{lowered}, into machine code through the use of general-purpose or domain-specific compilers, runtime calls to libraries, or a mix of multiple approaches. These transformations are designed to be sufficiently accurate---sometimes even exact---but may nevertheless change the derivatives. AD therefore does not necessarily differentiate ``what you implement'' but, rather, \emph{what is implemented at the level of abstraction at which AD is applied}, as we will illustrate in this section. It is generally the user's responsibility to ensure that the chosen abstraction's derivatives are reasonable approximations to those of the true function.
There are almost always levels of abstraction that are too low to permit this. At a hardware level, digital computers represent real or complex numbers using bit patterns, and functions as operations on these bit patterns. Viewed at this level, programs have discrete inputs and outputs and are therefore non-differentiable. Such a view is clearly unhelpful for computing derivatives; and AD \textit{usually} operates at a higher level of abstraction, where operations are assumed to represent functions with real (or complex) inputs and outputs~\cite{bolte2020mathematical}.

Most programs (including those written in relatively low-level languages such as C and Fortran) contain high-level operations that are implemented as sequences of lower-level operations within libraries or by the compiler, where the low-level and high-level semantics may differ in subtle ways. The mathematical function that we consider a program to implement may thus also depend on the compiler, math library, or target platform. Applying AD at a high abstraction level often avoids surprises arising from low-level implementation details but rests on the assumption that the high-level operations encountered at that level are evaluated exactly, a pitfall discussed in Section~\ref{sec:lowlevelaccuracy}. 

The abstraction that is most suitable for applying AD is not always explicitly represented in the source code. For example, a function involving multiple linear algebra operations can be implemented by using a sequence of loop nests within the same routine. In this case, the high-level linear algebra view that would often be most appropriate for AD~\cite{giles2008collected} exists only on paper (or in the mind of the programmer). Multiple levels of abstraction may exist in the same programming language, for example in the form of a C program calling a library function implemented in C, or in different programming languages, for example a Python script calling kernels written in C. Sometimes multiple abstraction levels are implemented manually, and sometimes they result from compiler lowering. One may therefore need to modify the implementation or build system to expose the appropriate abstraction in a language supported by a given AD tool. The remainder of this section provides examples for this pitfall.


\begin{figure}
    \centering
    \includegraphics[scale=0.7]{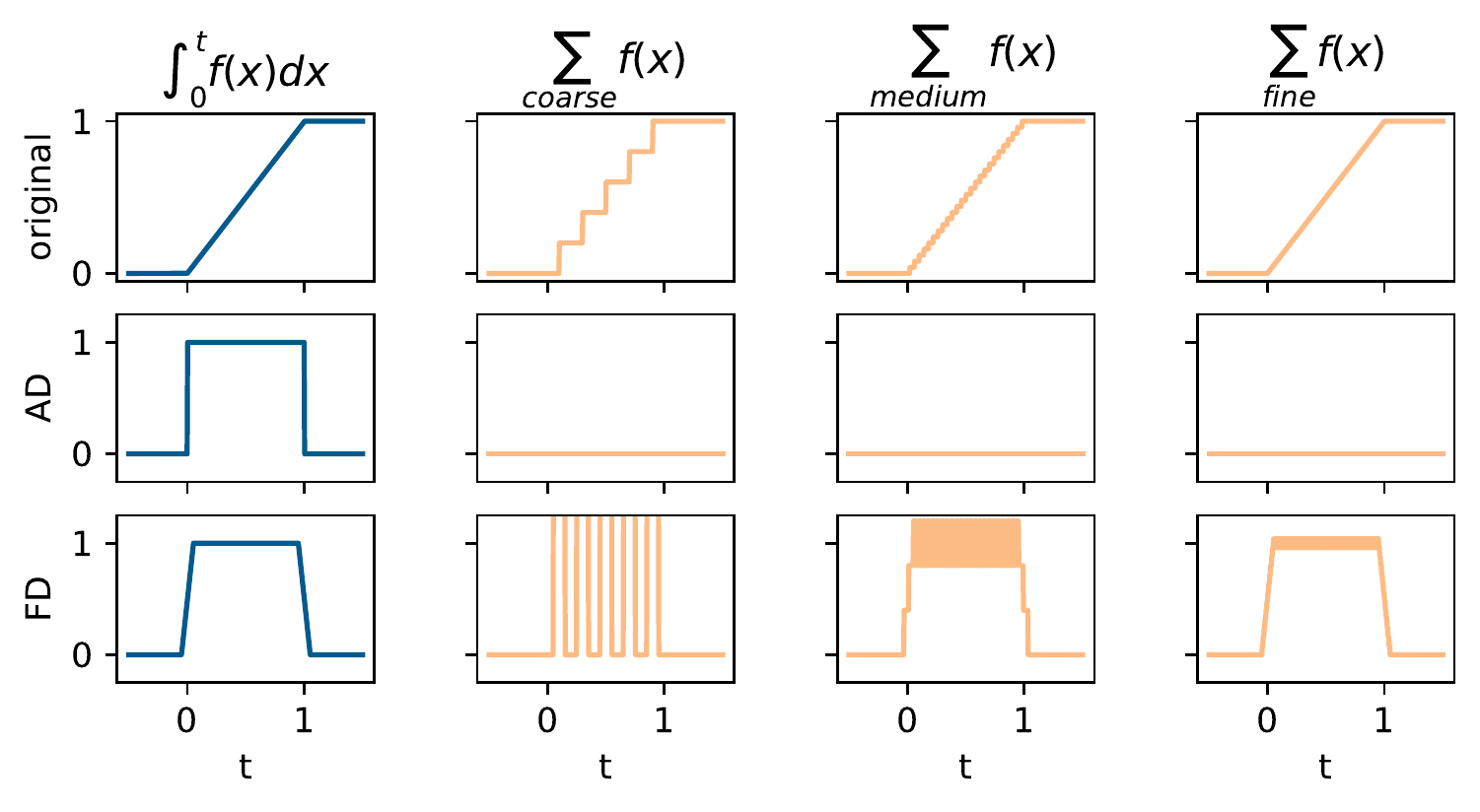}
    \caption{\textbf{Left column, blue:} Integration over a rectangular function yields the function in the top. Differentiating this integral by using AD or symbolic differentiation reconstructs the rectangular function, shown in the middle. Using finite differences (FD) yields a reasonable approximation, shown  at the bottom. \textbf{Other columns, orange:} If the integral is approximated by using quadrature rules, the function becomes a staircase with 0 derivatives almost everywhere. Increasing the resolution of the quadrature improves the accuracy of the primal but not that of the derivatives. Finite differences produce incorrect, but perhaps more intuitive, derivative results that converge to the rectangular function as the resolution is increased.}
    \label{fig:integral}
\end{figure}

\textbf{Explicit functions}
 are sometimes replaced with approximations whose derivatives differ from those of the approximated function. As an example, consider a call to the sine function $\sin(x)$ whose derivative involves a call to $\cos(x)$. At a lower abstraction level, however, the sine function might be implemented as a lookup table---possibly even containing bitwise-exact results for every floating-point input. The function is thus piecewise constant, with zero derivatives almost everywhere. Unexpected AD results caused by table lookups are known to occur in practice~\cite{ad-nuclear}. Other approximation techniques including bit hacks, range-reduction techniques, or bitwise reconstruction~\cite{863031,cordic} also have derivatives that are completely different from those of the approximated function.
A more benign example is the approximation of trigonometric functions using polynomials, whose derivatives are a polynomial of lower degree with plausible but less accurate results, especially for higher derivatives. Some programs approximate the same function by using different algorithms depending on the input value~\cite{gal1991accurate}, and differentiating through such an implementation yields useful derivatives for some but not all inputs.
AD practitioners may experience these problems with legacy code containing approximations or with math functions that are inlined by the build system before AD is applied. Approximations are not  used just to replace small math operations; they sometimes replace entire physical models, for example to incorporate empirical observations~\cite{ahmed2009search}, potentially leading to meaningless derivatives.

A related problem occurs for function approximations that are based on numerical quadrature~\cite{bangaru2021systematically}, Monte Carlo, or other methods where the output is essentially a (weighted) count of the number of points for which a certain condition holds. In the limit of infinitely many samples, the output may be differentiable with respect to certain problem parameters, but for any finite number of samples the count does not depend in a differentiable way on the parameters unless the abstraction level is raised~\cite{bangaru2021systematically}. See Figure~\ref{fig:integral} for an illustration.

\begin{figure}[t]
    \centering
    \includegraphics[scale=0.7]{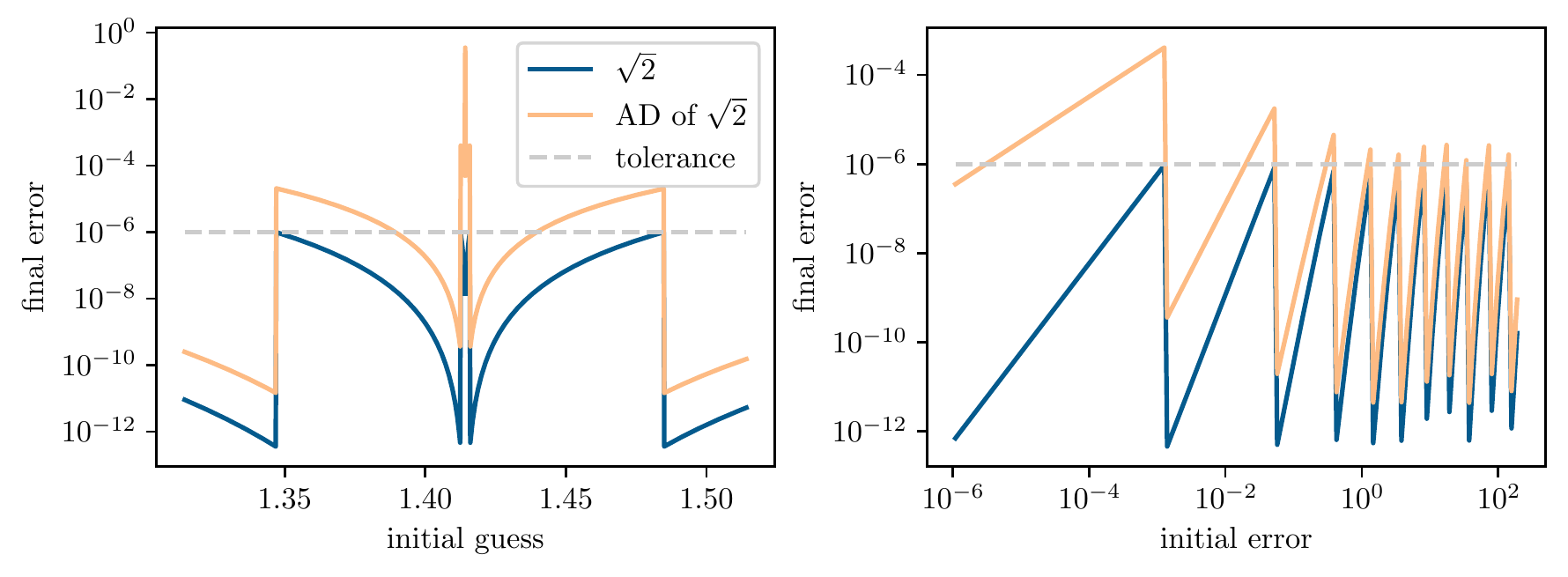}
    \caption{Square root finding with Heron's method for a set error tolerance of $10^{-6}$. When differentiating through the iterative method using AD, the stopping criterion of the primal is still used and stops the process before the derivatives reach the desired accuracy.}
    \label{fig:iterative}
\end{figure}

\textbf{Implicit functions} include the solution of systems of linear or nonlinear equations using direct or iterative processes.
Instead of differentiating through a linear solver implementation, it is usually more efficient and sometimes more accurate to formulate the derivative on a higher abstraction level in terms of a modified linear equation system~\cite{Griewank1993,racc,BARTHOLOMEWBIGGS2000171,MORE2014268}. When linear solvers are only incompletely converged, it can also be helpful to create an abstraction in which the initial guess and the final residual are considered as differentiable inputs and outputs~\cite{akbarzadeh2020consistent}.
Programs containing iterative fixed-point loops often require AD at a high abstraction level. For example, a user may want to obtain derivatives of the fixed-point location with respect to function parameters. The low-level implementation starts at some initial guess and refines the state until some accuracy or iteration count is reached. Especially when given a good initial guess, the process may terminate after few steps~\cite{Gilbert1992}. Accumulating the derivatives through this evaluation results in a high sensitivity to the initial state (which is actually $0$ for the true function $f$), and inaccurate derivatives with respect to function parameters. Figure~\ref{fig:iterative} shows an example in which the derivatives are orders of magnitude less accurate than the function evaluation, depending on the initial guess. Raising the level of abstraction allows us to view the process as an implicit function whose derivatives can be computed by using another iterative process with its own stopping criterion~\cite{BECK1994109,christianson1994reverse,christianson1998reverse}. This approach is implemented, for example, in the Tapenade AD tool~\cite{taftaf}, but it needs to be explicitly applied to individual loops by the user.


While we argue that the user is usually responsible for avoiding pitfalls related to the abstraction level, there are tools whose user interface and actual differentiation are at different levels of abstraction. For example, Enzyme~\cite{moses2020instead,moses2021reverse,10.5555/3571885.3571964} is called from the input language (e.g., C/C++, Julia), but the program is differentiated after lowering to LLVM-IR~\cite{LLVM}, which occasionally changes the derivatives. It is less clear who (the user, the language frontend, the AD tool) is responsible for ensuring correct derivative semantics in this case.


\subsection{Pitfall III: The chosen branch has unexpected derivatives}
\label{sec:branchdervs}

Chain rule differentiation of the operations encountered inside a branch may yield different derivatives from those encountered on other branches of the function. This can be problematic if program branches apply on closed subdomains---particularly for individual points. This problem has been discussed in~\cite[p.~343]{Griewank:2008:EDP:1455489} and~\cite{BECK1994119}.

Figure~\ref{fig:badbranch} shows a common pattern where the evaluation is simplified for certain input values. For example, BLAS functions~\cite{blas} often use a fast branch if one of the factors in a multiplication is $0$ or $1$, resulting in zero derivatives with respect to that factor. Since the branch is chosen only for individual points (a set of measure zero), the overall function and its derivatives should be unaffected. Unfortunately, AD will return the derivatives of the operations inside the branch, which could be any arbitrary value or could be correct, for example, for first-order derivatives but wrong for higher-order ones.
A similar problem occurs for \texttt{max} (or \texttt{min}) functions with inputs that contain more than one identical maximum (or minimum) value. The implementation may arbitrarily choose one of these values without impacting the function result, but the derivatives may change depending on the selected instance.

One might consider these problems less important since they  occur only at certain points, namely, on a measure zero space, and the probability of hitting such a ``bad spot'' might appear to be zero. However, previous work has pointed out that the probability is nonzero in practice because of finite precision~\cite[p. 336]{Griewank:2008:EDP:1455489} and because inputs are rarely randomly distributed and fast branches are often implemented precisely because they  speed up common cases.

\begin{figure}
    \centering
\noindent
\begin{tabular}{c|c|c}
\textbf{Definition} & \textbf{Code with fast path} & \textbf{Code with modified path} \\\hline
\begin{minipage}{12em}
\vspace{0.5em}Function:\\\vspace{-2em}
\begin{minted}[escapeinside=||]{python}
def f(x):
    return x

| |
\end{minted}
\end{minipage}
&
\begin{minipage}{12em}
\begin{minted}{python}
def g(x):
  if x == 0:
    return 0
  else:
    return x
\end{minted}
\end{minipage}
&
\begin{minipage}{15em}
\begin{minted}{python}
def h(x):
  if x == 0:
    return sin(x) # 0
  else:
    return x
\end{minted}
\end{minipage}
\\\hline
\begin{minipage}[c]{12em}
\vspace{0.5em}Derivative:\\\vspace{-2em}
\begin{minted}[escapeinside=||]{python}
def f_d(x, xd=1.0):
    return xd

| |
\end{minted}
\end{minipage}
&
\begin{minipage}{12em}
\begin{minted}{python}
def g_d(x, xd=1.0):
  if x == 0:
    return 0 # wrong!
  else:
    return xd
\end{minted}
\end{minipage}
&
\begin{minipage}{15em}
\begin{minted}{python}
def h_d(x, xd=1.0):
  if x == 0:
    return cos(x)*xd # xd
  else:
    return xd
\end{minted}
\end{minipage}
\\
\end{tabular} 
    \caption{The functions $f$, $g$, and $h$ return the exact same value for all inputs and implement the same mathematical function without any approximations. However, $g$ has a branch that, at $x=0$, will cause AD to compute different derivatives because the encountered operators do not depend on $x$ and thus have a $0$ derivative. The modified function $h$ illustrates that it is not the presence of a branch that causes a problem, but rather the inconsistency between derivatives in both branches. Further, the derivatives for $h$ change if a compiler optimization propagates the constant $x=0$ into the branch, effectively turning $h$ into $g$ before AD is applied.}
    \label{fig:badbranch}
\end{figure}

\subsection{Pitfall IV: The operators are not accurate enough}
\label{sec:lowlevelaccuracy}

\begin{figure}
    \centering
    \includegraphics[scale=0.7]{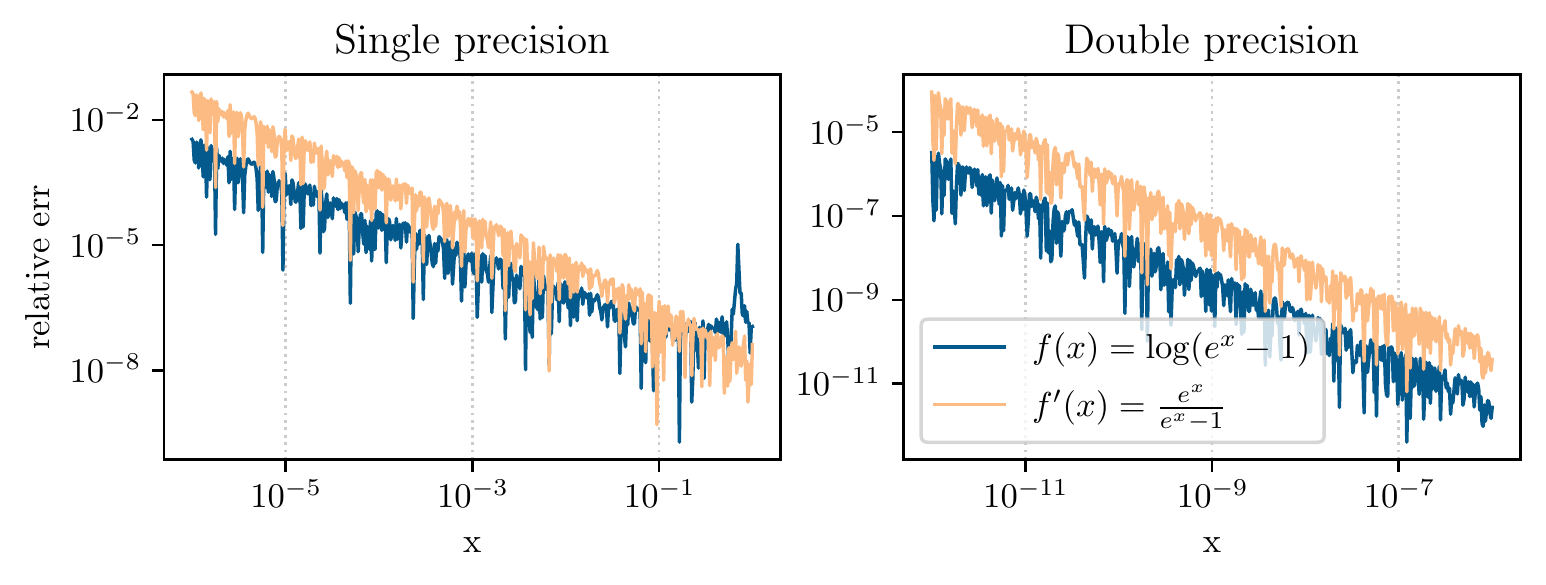}
    \caption{Example of a function whose derivative is orders of magnitude less accurate when evaluated with floating-point numbers, particularly for input values close to $0$.}
    \label{fig:accuracy_log}
\end{figure}

Sometimes, calculating the derivative of a function can incur more severe floating-point rounding errors than can calculating the function itself.
%
Consider the function $f(x) = \log(e^x - 1)$, whose derivative is given by $f'(x) = \frac{e^x}{e^x - 1}$.
Note that the derivative calculation involves dividing by $(e^x - 1)$. For small values of x, this is numerically unstable because it involves subtracting two nearly equal values, resulting in a division by a very small number with potentially large relative error.
Figure ~\ref{fig:accuracy_log} shows the errors introduced in the function and its derivative as $x$ approaches $0$.

Roundoff is not the only concern for this pitfall. Whenever AD is applied at a high level of abstraction, differentiation occurs with the assumption that operators at that abstraction level exactly compute some differentiable function whose derivatives can also be exactly computed by their corresponding derivative operators. This assumption is violated, for example, when iterative linear solvers are used that converge to a reasonably accurate value for the original function but break down, diverge, or converge more slowly for the derivative function. Commonly used iterative methods such as BiCG, BiCGSTAB, or GMRES with restarts do not have known performance guarantees and may break down in unpredictable ways~\cite{barrett1994templates}, or they  may succeed for a given matrix but fail for its transpose (as would be used for solving the adjoint system during reverse-mode AD), or they may be more likely to fail for certain right-hand side vectors that could be more likely to occur during derivative computations~\cite{huckelheim2017discrete}.


\section{Debugging Techniques}
\label{sec:solutions}

As alluded to in Section~\ref{sec:whattoexpect}, AD can be implemented in different modes that compute projections of Jacobian matrices or of their transpose, or higher-order derivative matrices.
Suppose an input program implements a function as a composition of two operations
$y = f(g(x))$, where $x$ and $y$ are real scalars or vectors of arbitrary size. Suppose that the Jacobian matrices of $f$, $g$ are $F$, $G$ respectively. For two given vectors $\dot{x}$ and $\bar{y}$, the forward and reverse modes of AD compute the projections
$\;\dot{y} = F G \dot{x}\;$ and $\;\bar{x} = \bar{y}G F,\;$
without explicitly forming $F$ or $G$.

\begin{figure}
    \centering
    \includegraphics[scale=0.7]{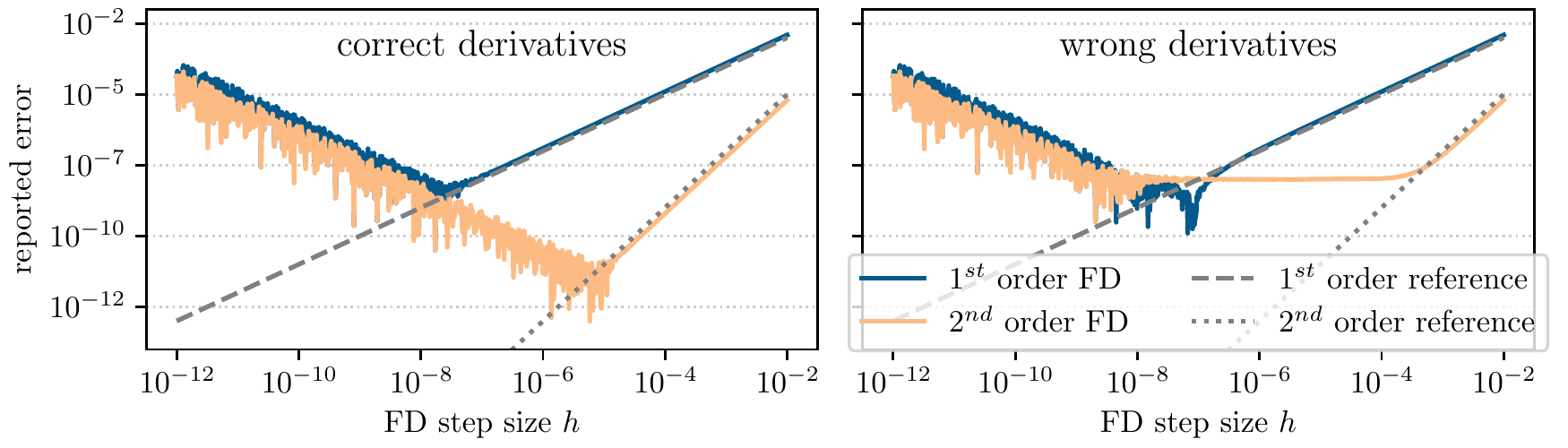}
    \caption{Comparing AD with first-order-accurate forward  and second-order-accurate central finite differences (FD). The step size trades roundoff (small $h$) for truncation errors (large $h$), with an implementation-dependent optimum. \textbf{Left:} Correct AD. \textbf{Right}: AD with $10^{-8}$ relative error, revealed by the second-order test over a range of $h$.}
    \label{fig:gradcheck}
\end{figure}

Finite differences (FD) can be used to obtain approximations for $\;\dot{y} = F G \dot{x}\;$ using $\dot{x}$ scaled with step size $h$ as input perturbation. This allows detection of AD tool bugs as well as many instances of pitfalls I, II and III, but  requires choosing a good step size and sometimes leaves users guessing whether a discrepancy is caused by floating-point or truncation errors or by AD problems~\cite{10.1145/356012.356013}.
Higher order formulae and multiple step sizes can be used for a more thorough test that determines the convergence order of the error as illustrated in Figure~\ref{fig:gradcheck}.

If a tool supports more than one mode or if multiple tools are available that support different modes, users can compare the results of these modes against each other.
One popular method for this purpose is the \emph{dot product test}. For the above example, we can compute the scalar quantity $\psi$ for some arbitrarily chosen $\dot{x}$ and $\bar{y}$ as follows: 
\begin{align*}
    \psi &= \bar{y}\dot{y}         && \text{definition}\\
         &= \bar{y}(G F\dot{x})    && \text{compute }\dot{y}\text{ with forward mode autodiff}\\
         &= (\bar{y}G) (F\dot{x})  && \text{intermediate step, using associativity}\\
         &= (\bar{y}G F)\dot{x}    && \text{compute }\bar{x}\text{ with reverse mode autodiff}\\
         &= \bar{x}\dot{x}         && \text{alternative definition}
\end{align*}
Computing $\psi$ with both forward- and reverse-mode AD allows us to detect inconsistencies between both modes. Decomposing the overall function and applying both modes to different function parts, as seen in the intermediate step, allows pinpointing the location of inconsistencies. While this is a useful debugging tool, pitfalls I, II and III would cause both modes to consistently produce the same wrong derivatives. 
To circumvent this problem, the dot product test can be modified by replacing forward mode AD with finite differences. This is implemented, for example, in the gradient testing routine in PyTorch~\cite{paszke2017automatic} using randomized values for $x$, $\dot{x}$ and $\bar{y}$, allowing direct debugging of reverse mode AD -- with the aforementioned caveats regarding step size and error tolerance. 

As with most tests, errors are only found if they affect the program for the tested inputs and program paths taken. In addition, the dot product test  checks a projection of the function's Jacobi matrix only onto the chosen vectors $\dot{x}, \bar{y}$ and may therefore miss error modes that are (almost) orthogonal to these vectors, especially if the test is implemented with a tolerance for roundoff errors.

\section{Conclusion}
We show in this paper that AD may produce wrong, or at least, surprising results. 
We believe that this need not be a reason to avoid AD. Many other methods---floating-point arithmetic, optimization algorithms, machine learning, to name a few---are useful despite occasional surprises. However, it is essential that users develop a mental model for the behavior of these techniques in order to avoid and detect issues. By expanding the list of  AD problems from previous work~\cite{Fischer1991SPi,BECK1994119,Griewank:2008:EDP:1455489,christianson1994reverse,bangaru2021systematically} and categorizing them into pitfalls, we hope to help users develop such a mental model.

At the same time, more research is needed to create languages and tools that implement AD in a more predictable way and make it easier to express programs in abstractions that have intuitive derivatives.
Recent work discusses AD semantics and provably correct AD but focuses on functional or domain-specific languages that are more restrictive than the languages used in most practical applications and that assume smoothness properties and real arithmetic instead of floating-point numbers~\cite{correctness-autodiff,provably-correct-higher-order,demystifying}.
There has also been recent progress on understanding the semantics of AD for neural networks that contain nonsmooth activation functions or that use finite precision arithmetic~\cite{NEURIPS2020_4aaa7617,lee2023correctness} or for probabilistic programs where users want to get derivatives of expectations with respect to input parameters~\cite{10.1145/3571198}. 

Programmers are often aware that they are responsible for the semantics of their program and that the language only guarantees semantics of individual constructs. We argue that a similar view of AD should be encouraged rather than making sweeping claims about AD ``differentiating programs.'' The quote from~\cite{naumann2012art} remains true a decade later:
\emph{The application of AD to computer programs still deserves to be called an “art.”}

\section*{Acknowledgments}
\small{
We thank Christian Bischof, Valentin Churavy, and Jesse Michel for insightful comments and for bringing pitfall examples to our attention. We are grateful for discussions with the late Andreas Griewank, who  greatly influenced our work.
This work was supported by the Applied Mathematics activity within the U.S. Department of Energy, Office of Science, Advanced Scientific Computing Research Program, under contract number DE-AC02-06CH11357, Lawrence Livermore National Laboratory under Contract DE-AC52-07NA27344, LANL grant 531711, and by the US Air Force Research Laboratory and the United States Air Force Artificial Intelligence Accelerator under Cooperative Agreement Number FA8750-19-2-1000. The views and conclusions contained in this document are those of the authors and should not be interpreted as representing the official policies, either expressed or implied, of the United States Air Force or the U.S. Government. The U.S. Government is authorized to reproduce and distribute reprints for Government purposes notwithstanding any copyright notation herein.
}


\printendnotes

\bibliography{pitfalls}



\end{document}